\documentclass[12pt,amsfonts,amscd]{amsart}
\usepackage{amssymb}
\newtheorem{theorem}{Theorem}[section]
\newtheorem{lemma}[theorem]{Lemma}
\newtheorem{prop}[theorem]{Proposition}
\newtheorem{corollary}[theorem]{Corollary}

\theoremstyle{definition}
\newtheorem{definition}[theorem]{Definition}
\newtheorem{example}[theorem]{Example}

\theoremstyle{remark}

\numberwithin{equation}{section}

\newcommand{\DD}{{\mathbb D}}

\newcommand{\RR}{{\mathbb R}}
\newcommand{\CC}{{\mathbb C}}

\DeclareMathOperator{\psh}{PSH}

\renewcommand{\phi}{\varphi}

\begin{document}
\title
{Graphs with multiple sheeted pluripolar hulls}

\author{Evgeny Poletsky}

\address{Department of Mathematics,  215 Carnegie Hall,
Syracuse University,  Syracuse, NY 13244} \email{eapolets@syr.edu}

\author{Jan Wiegerinck}

\address{Korteweg--de Vries institute for mathematics,
University of Amsterdam, Plantage Muidergracht 24, 1018 TV,
Amsterdam, The Netherlands} \email{janwieg@science.uva.nl}






\begin{abstract} In this paper we study the pluripolar hulls of
analytic sets. In particular, we show that hulls of graphs of analytic
functions can be multiple sheeted and sheets can be separated by a
set of zero dimension.
\end{abstract}
\thanks{The first author was partially supported by an NSF Grant.}
\maketitle

\section{Introduction}\label{S:in}
One of the oldest interesting topics in complex analysis is
the problem of analytic extensions: find the maximal analytic
object containing a given one. For example, if $f$ is an
analytic function we are looking for its analytic continuation and
if $A$ is an irreducible analytic set we try to find another one of the same dimension
containing $A$.
\par The counterpart of analytic
extension in pluripotential theory is the so-called {\it pluripolar hull.} There are two types of
{\em pluripolar hulls} of a set $A$ in a domain $D\subset\CC^n$.
Let $\psh(D)$ be the set of all plurisubharmonic functions on $D$
and $\psh_0(D)$ the set of all negative functions from
$\psh(D)$. Define
$$A^*_D=\{z\in D: \forall h\in \psh(D)\, h|_A=
-\infty \Rightarrow h(z)=-\infty\}$$ and
$$A^-_D=\{z\in D: \forall h\in \psh_0(D)\, h|_A=
-\infty \Rightarrow h(z)=-\infty\}.$$
 For example, if $A$ is an
analytic set in  pseudoconvex domain $D$, then every point of
$A$ has a neighborhood $V$ where $A\cap V=\{h_1=\dots=h_k=0\}$ and
the functions $h_k$ are holomorphic on this neighborhood. Hence
$A\cap V=\{\log\max\{|h_1|,\dots,|h_k|\}=-\infty\}$. In fact, $A^*_D=A$,
 because by \cite[Cor.
1]{C} there even exists $v\in \psh(D)$ such that
$A=\{v=-\infty\}$. 

If such a $v\in\psh(D)$ exists for $A$, we call $A$ {\it
pluricomplete} in $D$.  In general, an
analytic extension of $A$ is contained in $A^*_D$.
\par In the case when $A=\Gamma_f$ is the graph of an analytic
function $f$ it was boldly conjectured in \cite{LMP} that the
closure of the analytic extension of $A$ coincides with $A^*_D$.
However, A. Edigarian and the second author found in \cite{EW1} an
analytic function $f$ on the unit disk that does not extend
analytically while the pluripolar hull of its graph is a graph of a function
defined on almost the whole plane, cf.~\cite{siciak}. This example raised the question: what
are pluripolar hulls or, better to say, extensions of analytically
non-extendible analytic sets?
\par The pluripolar hull of the graph $\Gamma_f$ of a holomorphic
function $f(z)$ on a domain $D$ may well be multi-sheeted over
$D$. The principal value of $\sqrt z$ on $\{\Re z>0\}$ provides
the easiest example. Only recently Zwonek \cite{Zwonek}, and,
independently, Edlund \& J\"oricke \cite{Edl-Jo} gave examples of
holomorphic functions $f$ on their {\em  domain of existence} $D$
with the property that the pluripolar hull $(\Gamma_f)^*_{\CC^2}$
is multi-sheeted over (parts of) $D$.
\par In these examples sheets can be separated only by a cut whose
projection on $z$-plane has the dimension 1. As we show in Section
\ref{S:pe} this is an intrinsic property of analytic extensions.
\par In the present note we show that there is a Cantor type set
$K$ and a holomorphic function $f(z)$ on $D=\CC\setminus K$ such
that $(\Gamma_f)^*_{\CC^2}$ is 2-sheeted over $D$. So for
pluripolar extensions sheets can be separated by a $0$-dimensional
cut. As a by-product we obtain an example of a uniformly
convergent sequence of holomorphic functions such that their
pluripolar hulls do {\em not} converge to the pluripolar hull of
the limit.
\par The set $K$ should be sufficiently fat. Edigarian and
the second author showed that if $D=\CC\setminus K$, with $K$ a
polar compact set in $\CC$, and if $f$  is not extendible over
$K$, then $(\Gamma_f)^*_{\CC^2}\cap D\times \CC=\Gamma_f$; see
\cite{EW2}, and \cite{EW3} for the fact that also over $K$ the
hull is at most single sheeted.
\section{Pluripolar extensions}\label{S:pe}
\par Let $E$ be a closed set in a pseudoconvex domain
$D\subset{\CC}^n$. If $A\subset D\setminus E$ then, in general,
$A^*_{D\setminus E}$ is a proper subset of $A^*_D\setminus E$.
However, as the following statement shows, these sets coincide
when $E$ is pluripolar.
\begin{prop}\label{p:ps} If $E$ is a closed pluripolar set in a pseudoconvex
domain $D\subset{\CC}^n$ and $A\subset D\setminus E$ then,
$A^*_D\setminus E= A^*_{D\setminus E}$.
\end{prop}
\begin{proof} Let $\{D_j\}$ be an increasing sequence of
relatively compact subdomains with $\cup_jD_j=D$. By \cite[Thm.
2.4]{LP} $A^*_D=\cup_j(A\cap D_j)^-_{D_j}$. If
$u\in\psh_0(D_j\setminus E)$, then $u$ extends as a negative
plurisubharmonic function to $D_j$ (see \cite[Thm. 2.9.22]{K}).
Therefore, $(A\cap D_j)^-_{D_j}\setminus E=(A\cap
D_j)^-_{D_j\setminus E}$.  Since $(A\cap D_j)^-_{D_j\setminus
E}\subset A^*_{D\setminus E}$, we see that $A^*_D\setminus
E\subset A^*_{D\setminus E}$ and, consequently, $A^*_D\setminus E=
A^*_{D\setminus E}$.
\end{proof}
\par The proposition below describes the situation when $E$ is a
closed set in a pseudoconvex domain $D\subset{\CC}^n$, $A\subset
D\setminus E$ is an analytic set and $A^*_D$ is also analytic.
\begin{prop}\label{p:dim} Suppose that $E$ is a closed set in
a pseudoconvex domain $D\subset{\CC}^n$ and $A\subset D\setminus
E$ is an analytic set. If the set $A^*_D$ is analytic then
every irreducible component of $A^*_D$ contains a component of $A$ of the
same dimension.
\end{prop}
\begin{proof} Let $X$ be an irreducible component of $A^*_D$. We
represent $A^*_D$ as $X\cup Y$, where $Y$ is another analytic set
in $D$ and $\dim X\cap Y<\dim X$. As we indicated in Section
\ref{S:in} the set $Y$ is pluricomplete and $Y^*_D=Y$. It is easy
to check that if sets $F,G\subset D$, then $(F\cup
G)^*_D=F^*_D\cup G^*_D$. So if $B=A\cap(X\setminus Y)$ then
$X\setminus Y\subset B^*_D$.
\par Suppose that $\dim X>\dim B$ and let $R$ be the set of
regular points of $X$. The set of singular points of $X$ is
analytic and, consequently, pluricomplete. By the argument above
$X\setminus Y$ belongs to the pluripolar hull of the set $B'=B\cap
R$.
\par We may assume that $0\in R$ and let $T$ be the tangent plane
to $R$ at 0.  If $p$ is a projection of $\CC^n$ on $T$, then the
set $p(B')$ is pluripolar in $T$ and, consequently, there is a
plurisubharmonic function $u$ on $T$ equal to $-\infty$ on
$p(B')$. The set $p(R)$ has a non-empty interior in $T$ and,
therefore, there is a point $z_0\in R$ such that
$u(p(z_0))\ne-\infty$. Then the function $v=u\circ p$ is
plurisubharmonic on $\CC^n$, equal to $-\infty$ on $B'$ and
$v(z_0)>-\infty$. Thus $R$ does not belong to the pluripolar hull
of the set $B'$. This contradiction proves the proposition.
\end{proof}
\par Suppose that $A$ is a pluricomplete analytic set of pure
dimension $m$ in $D\setminus E$. If $A^*_D$ is an analytic set in
$D$ and $A$ is a proper subset of $A^*_D\setminus E$, then the set
$E\cap A^*_D$ cuts $A^*_D$ into several pieces and, therefore, its
topological dimension must be at least $2m-1$.
\par  For example, let
$D=\{(z,w)\in{\CC}^2\}$ and $E=\{\Im z=0, \Re z\ge0\}$. Take a
branch $w=f(z)$ of the function $w=\sqrt{z}$ over $\CC \setminus
E$ and let $A=\{(z,f(z)):\, z\in{\mathbb C}\}$. The pluripolar
hull of $A$ in ${\CC}^2\setminus E$ is $A$ because the function
$\log|f(z)-z|$ is equal to $-\infty$ exactly on $A$. But
$A^*_D=\{(z,w):\,z=w^2\}$. In this example $A^*_D$ is an analytic
set and the set $A^*_D\cap E$ is the real curve $\{(x^2,x):\ x\in
\RR\}$ which projects 2 to 1 except at 0 and its projection has
dimension 1.

\par As the following statement shows this is the minimal possible
dimension.
\begin{prop} Let $E$ be a closed set in a pseudoconvex domain
$D\subset{\CC}^n$ and let $A$ be an irreducible analytic set of
dimension $m$ in $D\setminus E$ such that $A^*_D$ is also
analytic. If $p$ is a projection of $\CC^n$ onto
$\CC^m\subset\CC^n$ such that $p(A)$ has a non-empty interior in
$\CC^m$ and the topological dimension of $p(E)$ is less than
$2m-1$, then $A^*_D\setminus E=A$.
\end{prop}
\begin{proof} By Proposition \ref{p:dim} every irreducible
component of $A^*_D$ contains a component of $A$ of the same
dimension. Thus $A^*_D$ is also irreducible and has dimension $m$.
We denote by $X$ the set of regular points of $A^*_D$ and let $X'$
be the subset of $X$ where the restriction of the projection $p$
to $X$ has maximal rank $m$. Since $p(A)$ has non-empty
interior in $\CC^m$, $x'$ is non-empty and relatively open in $X$. The
set $X\setminus X'$ is analytic and, therefore, has dimension at
most $m-1$. Hence the set $X'$ is connected and the set $A\cap X'$
is not empty.
\par Suppose that the topological dimension of $p(E)$ is smaller
than $2m-1$. Choose points $z_0\in A\cap X'$ and $z_1\in X'$ such
that $p(z_0)$ and $p(z_1)$ do not belong to $p(E)$. We can connect
these points by a real analytic curve $\gamma$ in $X'$. The upper
bound on the dimension of $p(E)$ implies that we can slightly shift
$\gamma$ so that $p(\gamma)$ does not meet $p(E)$. Since
$E$ is closed there is a relatively open neighborhood of $\gamma$
in $X'$ which does not contain points of $E$. But $A$ is a
relatively open analytic subset of $X'\setminus E$. Hence
$\gamma\subset A$ and $z_1\in A$.
\par If $p(z_1)\in p(E)$ but $z_1\not\in E$, then we can take a
neighborhood of $z_1$ in $X'$ where points $z$ with $p(z)\in E$
form a set with empty interior. Since other points are in $A$, the
whole neighborhood is there also. Thus $X'\setminus E\subset A$.
Since $A$ is closed in $A^*_D\setminus E$, $A^*_D\setminus E=A$.
\end{proof}
\par When $A$ is an analytic set in $D\setminus E$ we denote by
$A_E$ the intersection of the closure of $A$ in $D$ with $E$. If
$A^*_D\setminus E\ne A$ we will say that $A$ has a {\it
non-trivial pluripolar extension } through $E$ in $D$.
\par The following theorem lists some limitations on the
set $A_E$ when a non-trivial pluripolar extension takes place.
Following \cite{Ch} we call a set $G$ in a domain $Y\subset{\mathbb
C}^p$ {\it locally removable } if $G$ is closed and for every open
set $V$ in $Y$ every bounded holomorphic function $f$ on
$V\setminus G$ extends holomorphically to $V$.
\begin{theorem}\label{t:le} Suppose that $D$ is a pseudoconvex domain
in ${\mathbb C}^n$, $E$ is a closed set in $D$ and $A$ is an analytic
set of pure dimension $m$ in $D\setminus E$ with a non-trivial
pluripolar extension through $E$ in $D$. Then the
$(2m-1)$-Hausdorff measure of the set $A_E$ is not equal to zero
and if, additionally, $p:\,{\mathbb C}^n\to{\mathbb C}^m$ is a
projection such that the restriction $p|_{A_E}$ is proper and
$A\cap p^{-1}(z)$ is empty for all $z\in p(A_E)$, then $p(A_E)$ is
not locally removable in ${\mathbb C}^m$.
\end{theorem}
\begin{proof} The set $A$ is analytic in $D\setminus A_E$. If the
$(2m-1)$-Hausdorff measure of $E$ is  zero, then by Shiffman's
theorem (see \cite[4.4]{Ch}) the closure $\bar A$ of $A$ in $D$ is
an analytic set in $D$. Since the domain $D$ is pseudoconvex there
is a holomorphic function $f$ on $D$ such that $\bar A=\{f=0\}$.
Thus $A^*_D=\bar A$  and this contradicts the assumption that the
extension is non-trivial.
\par In the second case if $p(A_E)$ is locally removable
in ${\mathbb C}^m$, then by the proposition in \cite[18.1]{Ch} the
closure $\bar A$ of $A$ is an analytic set in $D$ as before and
the same argument leads to a contradiction.
\end{proof}
\par In our main example $n=2$, $m=2$ and $m=1$. In this case
Theorem \ref{t:le} can be reformulated as follows:
\begin{corollary} If in the assumptions of Theorem \ref{t:le} $n=2$
and $m=1$, then the first Hausdorff measure of $A_E$ is not zero
and, under additional assumptions, the first Hausdorff measure of
$p(A_E)$ are not zeros.
\end{corollary}
\section{A holomorphic function on the complement of a Cantor type set with
2-sheeted hull}
\begin{definition}
A {\em Cantor type set} $K$ will be  a compact perfect
subset of $\RR$ with empty interior.
\end{definition}
It is a well known fact from elementary point set topology that
such a $K$ is homeomorphic with Cantor's middle third set. It is
of the form $[a_0,b_0]\setminus (\cup_{j=1}^\infty  I_j)$ where
$I_j$ are open intervals in $[a_0,b_0]$,
$\overline{I_j}\cap\overline{I_k}=\emptyset$ if $j\ne k$ and $
\cup_{j=1}^\infty  I_j$ is dense in $[a,b]$. We can assume that
the length of $I_j=(a_j,b_j)$ decreases with $j$.

It is useful to enumerate the set $\{a_j,b_j,\ j=0,\ldots, n\}$ as
$\{\alpha_{jn},\beta_{jn}\}$ so that $\alpha_{0n}=a_0$,
$\alpha_{jn}<\beta_{jn}<\alpha_{j+1,n}$ and $\beta_{nn}=b_0$. Note
that $[a_0,b_0]\setminus\cup_{j=1}^nI_j=
\cup_{j=0}^n[\alpha_{jn},\beta_{jn}]$ and that $I_m\cap
[\alpha_{jn},\beta_{jn}]\ne\emptyset$ implies that
 $I_m\subset [\alpha_{jn},\beta_{jn}]\ne\emptyset$.

\par Let $$g_n(z)=\frac{z-a_0}{z-b_0}\frac{z-b_1}{z-a_1}\cdots
\frac{z-b_n}{z-a_n}.$$
Then
$$g_n(z)=\frac{z-\alpha_{0n}}{z-\beta_{0n}}
\frac{z-\alpha_{1n}}{z-\beta_{1n}}\cdots
\frac{z-\alpha_{nn}}{z-\beta_{nn}}.
$$
Each fraction $\frac{z-\alpha_{jn}}{z-\beta_{jn}}$,
$j=0,1,\dots,n$, has a holomorphic branch
$\sqrt{\frac{z-\alpha_{jn}}{z-\beta_{jn}}}$ of its square root
outside $[\alpha_{jn},\beta_{jn}]$ that equals  1 at infinity. Let
\begin{equation}
\label{eq1} f_n(z)=\sqrt{\frac{z-\alpha_{0n}}{z-\beta_{0n}}}
\sqrt{\frac{z-\alpha_{1n}}{z-\beta_{1n}}}\cdots
\sqrt{\frac{z-\alpha_{nn}}{z-\beta_{nn}}}=\sqrt{g_n(z)}.
\end{equation}
Then $f_n(\infty)=1$, $f_n^2=g_n$ and $f_n$ is holomorphic on
$G_n=(\CC\setminus[a_0,b_0])\cup_{j=1}^nI_j$.

The maximal analytic extension of $f_n$ is a branched two sheeted
cover $X_n=\{(z,w): w^2=g_n(z)\}$ of $\CC$ that branches over
$\{a_j,b_j:\ j=0,1,\ldots, n\}$. The pluripolar hull
$(\Gamma_{f_n})^*_{\CC^2}$ equals $X_n$.
\begin{lemma}\label{lem} Keeping the notation as above, the sequence $\{g_n\}$
 converges normally to an analytic function $g$ on $\CC\setminus
K$. Moreover, the function $g$ extends analytically over a point
$x\in K$ if and only if for some $\alpha<x<\beta$ the length of
the set $K\cap(\alpha,\beta)$ is zero.
\end{lemma}
\begin{proof} Let $L$ be a compact set in $\CC\setminus K$. Let us show that
\begin{equation}
\label{eq2}\frac{z-a_0}{z-b_0}\prod_{j=1}^\infty\frac{z-b_j}{z-a_j}=
\lim_{n\to\infty}g_n(z)
\end{equation}
is uniformly convergent on $L$. There exists $n_0$ such that
$L\subset G_n$ for $n>n_0$ and moreover, for some $\delta>0$
$$L\subset\{z: |z-a_j|>\delta, j=0,1,\dots\}.$$
Hence, for  $z\in L$,
$$\left|\frac{z-b_j}{z-a_j}-1\right|=
\left|\frac{b_j-a_j}{z-a_j}\right|\le\frac{b_j-a_j}{\delta} .$$
Since $\sum(b_j-a_j)$ is finite, the product in \eqref{eq2}
converges uniformly on $L$ to a function $g$ that is holomorphic
on $\CC\setminus K$.
\par Suppose that the function $g$ extends analytically over a point $x\in
K$ so that $g$ is analytic on $(\CC\setminus
K)\cup(\alpha,\beta)$. We may assume that $\alpha\in I_k$, $b_k\le
x$, and $\beta\in I_m$, $a_m\ge x$ and
$$g_{1n}(z)=\prod\frac{z-\alpha_{jn}}{z-\beta_{jn}},$$
where the product runs over all $j$ such that either
$\beta_{jn}<\alpha$ or $\alpha_{jn}>\beta$. Let
$$g_{2n}(z)=\prod\frac{z-\alpha_{jn}}{z-\beta_{jn}},$$ where the
product runs over all $j$ such that $\alpha<\alpha_{jn}$ and
$\beta_{jn}<\beta$. Then $g_n=g_{1n}g_{2n}$ and by the argument
above the sequences $\{g_{1n}\}$ and $\{g_{2n}\}$ converge
uniformly on compacta on $\CC\setminus(K\setminus(\alpha,\beta))$
and $\CC\setminus(K\cap(\alpha,\beta))$ respectively. We denote
their respective limits by $g_1$ and $g_2$.
\par The derivative
$g_{2n}(\infty)=\sum(\beta_{jn}-\alpha_{jn})=l_n$, where the sum
runs over all $j$ such that $\alpha<\alpha_{jn}$ and
$\beta_{jn}<\beta$. Thus $g_{2n}(\infty)$ is equal to the length
$l_n$ of the intervals $(\alpha_{jn},\beta_{jn})$ lying in
$(\alpha,\beta)$ and $g'_2(\infty)$ is the length of the set
$K\cap(\alpha,\beta)$. If this length is positive, then the
function $g_2$ is not constant and, therefore, does not extend to
$K\cap(\alpha,\beta)$.
\par If this length is 0 then for $z\in\CC$ such that $|z-y|\ge1$
for all $y\in(\alpha,\beta)$ we have
$$\left|1-g_{2n}(z)\right|=
\left|1-\prod\left(1+\frac{\beta_{jn}-\alpha_{jn}}{z-\beta_{jn}}\right)\right|
\le e^{l_n}-1.$$ Hence the sequence $\{g_{2n}\}$ converges to 1
near $\infty$, $g_2\equiv1$ and $g$ extends analytically over
$(\alpha,\beta)$.
\end{proof}
\begin{lemma}\label{lem1}
If $f=f_K$ and the length of $K$ is positive, then the union
$\bar\Gamma_f\cup\bar\Gamma_{-f}$ of the closures of the graphs of
$f$ and $-f$ is not an analytic set.
\end{lemma}
\begin{proof}
\par If $A=\bar\Gamma_f\cup\bar\Gamma_{-f}$ is an analytic set, then
there is a holomorphic function $h=h(z,w)$ on $\CC^2$ such that
$h\equiv0$ on $A$. We have $A=\Gamma_f\cup\Gamma_{-f}\cup E$,
where $E\subset K\times\CC$.
\par Let us show that for every $z_0\in K$ the analytic set
$E_{z_0}=\{w:\,(z_0,w)\in E\}=A\cap\{z_0\}\times\CC$ consists of
at most two points. If it contains three points, then at least two
of them belong to, say, $\bar\Gamma_f$. Thus there are sequences
$\{z_j\}$ and $\{z'_j\}$ converging to $z_0$ such that the
sequences $\{f(z_j)\}$ and $\{f(z'_j)\}$ have distinct limits.
Connecting each $z_j$ and $z'_j$ by small curves in $\CC\setminus
K$ and looking at their limits we see that the cluster set of $f$
at $z_0$ contains a continuum. Hence, $E_{z_0}=\CC$ and
$h(z_0,w)\equiv0$.
\par From the Taylor expansion of $h$ we immediately derive that
$h(z,w)=(z-z_0)^nh_1(z_0,w)$, where $h_1$ is holomorphic on
$\CC^2$ and $h_1(z_0,w)\not\equiv0$. But for every point $w\in\CC$
there is a sequence of $z_j$ converging to $z_0$ such that, say,
$f(z_j)$ converges to $w$. Since $h(z_j,f(z_j))=0$ we see that
$h_1(z_0,w)=0$. This contradiction shows that $E_{z_0}$ has at
most two points and the intersection of $\bar\Gamma_f$ or
$\bar\Gamma_{-f}$ with $E$ consists of at most one point.
\par It follows that $f$ extends continuously to $K$. Since $K$
lies on the real line $f$ extends holomorphically to $\CC$ but
this impossible by Lemma \ref{lem}.
\end{proof}
\begin{example}
If the set $K$ has Lebesgue-measure 0, then
$$\lim_{n\to\infty} g_n(z)=1,$$
uniformly on any compact set $L$ not meeting $K$. It follows that
$f\equiv 1$ and $(\Gamma_f)^*_{\CC^2}=\{(z,1)\}$. But the
Hausdorff limit of the sets $X_n$ over $D$ equals $\{(z,w):\, w=1
\text{\ or \ } -1\}$.
\end{example}
\par We will need the next lemma whose proof is similar to the
proof of Theorem 2.1 in \cite{EW2}.
\begin{lemma}\label{fa}
Let $f$ be a holomorphic function on a domain $V\subset\CC^n$
containing a closed ball $B$ and let $\{r_n\}_n$ be a sequence of
rational functions of degree $n$ with poles outside $V$ and such
that the sup-norm $\|f-r_n\|_B^{1/n}\to0$ as $n\to\infty$. Then
there is a plurisubharmonic function $v$ on $\CC^{n+1}$ such that
$\{v=-\infty\}\cap(V\times\CC)=\Gamma_f$. Thus,
$(\Gamma_f)^*_{\CC^{n+1}}\cap(V\times\CC)=\Gamma_f$.
\end{lemma}
\begin{proof} The functions $r_n$ are ratios of polynomials $p_n$
and $q_n$ of degree $n$. We may assume that $B$ is the closed unit
ball centered at the origin and $\|q_n\|_B=1$. Then
$|q_n(z)|\le\max\{1,|z|^n\}$, $\|p_n\|_B$ does not exceed some
constant $C$ and $|p_n(z)|\le C\max\{1,|z|^n\}$.
\par Consider the plurisubharmonic functions
$$u_n(z,w)=\frac1n\log|q_n(z)w-p_n(z)|$$ on $\CC^{n+1}$. From the
estimates on $p_n$ and $q_n$ there is a constant $C_1$ such that
$u_n(z,w)\le2\log|z|+\log|w|+C_1$ when $|z|\ge1$ and
$u_n(z,w)\le\log|w|+C_1$ when $|z|\le1$.
\par We take $z_n\in B$ such that $q_n(z_n)=a_n$,
$|a_n|=1$ and let $w_n=(p_n(z_n)+1)/a_n$. Then $|w_n|\le C+1$ and
$u_n(z_n,w_n)=0$. If $B_n$ is a ball in $\CC^{n+1}$ centered at
$(z_n,w_n)$ and of radius $r_n=C+5$, then $B_n$ contains the unit
ball $B'$ centered at the origin and
$$\int\limits_{B_n}u_n\,dV\ge cu_n(z_n,w_n)=0.$$
It is immediate from the upper estimates on $u_n$ that there is
a constant $C_2$ such that
$$\int\limits_{B'}u_n\,dV\ge C_2,$$

\par By our assumption there is a sequence $\{d_n\}$ converging to
$\infty$ such that
$$u_n(z,f_n(z))=
\frac1n\log|q_n(z)|+
\frac1n\log\left|f(z)-\frac{p_n(z)}{q_n(z)}\right|\le
-d_n$$when $|z|\le1$.
\par Let us take a sequence $\{c_n\}$ of positive reals such that
$\sum c_n=1$ while $\sum c_nd_n=\infty$. Let
$$v(z,w)=\sum_{n=1}^\infty c_n\max\{u_n(z,w),d_n\}.$$
Since
$$\int\limits_Bv\,dV\ge C_2,$$
$v\not\equiv-\infty$ and, therefore, is a plurisubharmonic
function on $\CC^{n+1}$. Clearly, $v(z,f(z))=-\infty$ when
$|z|\le1$. Therefore, $v=-\infty$ on $\Gamma_f$.
\par The zeros of the polynomials $q_n$ are in $\CC^n\setminus V$.
Hence the functions $h_n(z)=\frac1n\log|q_n(z)|$ are harmonic on
$V$ and uniformly bounded above on compacta. So if $z\in V$ and
$\liminf h_n(z)=-\infty$, then there is a subsequence $\{n_k\}$
such that $\lim h_{n_k}(z)=-\infty$. Therefore, functions
$h_{n_k}$ converge to $-\infty$ uniformly on compacta in $V$. But
$h_{n_k}(z_{n_k})=0$ and this contradiction tells us that $\liminf
h_n(z)>-\infty$. So if $w\ne f(z)$, then
$$v(z,w)\ge \sum_{n=1}^\infty c_nh_n(z)+
\sum_{n=1}^\infty
\frac{c_n}n\log\left|w-\frac{p_n(z)}{q_n(z)}\right|>-\infty.$$
Hence $\{v=-\infty\}\cap(V\times\CC)=\Gamma_f$ and
$(\Gamma_f)^*_{\CC^{n+1}}\cap(V\times\CC)=\Gamma_f$.
\end{proof}
\par Now we can present our main example.
\begin{theorem} \label{ex1} There exists a Cantor type set $K$
obtained by deleting intervals $I_i=(a_i,b_i)$ from $[-1,1]$, such
that the function $f=f_K$ given by Lemma \ref{lem1} has the
following property: $$(\Gamma_f)^*_{\CC^2}\cap(\CC\setminus
K)\times\CC=\Gamma_f\cup\Gamma_{-f}.$$
\end{theorem}

\begin{proof} We will construct $K$ by deleting a sequence of
open intervals $(a_i,b_i)$ from the interval $[-1, 1]$. For
convenience, set $a_0=1$, $b_0=-1$. In order to choose the
intervals appropriately, we have to construct certain subdomains
$D_n$ in the open unit disk $\DD$ in the process. The domains
$D_n$ will contain the set $\DD\cap \{|\Im z|>1/2\}$. Thus the
closed discs $S=\{|z+3i/4|\le 1/8\}$ and $X=\{|z-3i/4|\le 1/8\}$
will be contained in $D_n$.
\par For a compact set $F$ in a domain $D\subset\CC$ let
$$\omega(z,F,D)=
-\sup\{h(z): h\in \psh_0(D), \limsup_{w\to K}h(w)\le-1\}$$ be the
harmonic measure of $F$ in $D$. Set $D_0=\DD$ and observe that
$\omega(z,S,D_0)>c_0$ for some positive $c_0$.
\par Let $\{c_n\}$ be a sequence of positive real numbers
converging to $\infty$. Suppose that the intervals $I_1,\ldots,
I_n$ have been chosen. We take as $a_{n+1}$ the midpoint of the
largest interval in their complement. Next take $b_{n+1}>a_{n+1}$
so small that the interval $[a_{n+1},b_{n+1}]$ does not intersect
the intervals $I_1,\ldots, I_n$,
$d_{n+1}=b_{n+1}-a_{n+1}<4^{-(n+1)c_{n+1}}$ and, moreover,
$$\omega(z, S, D_{n+1})>c_0, \quad z\in X.$$
Here we define $D_{n+1}=D_n\setminus \DD(a_{n+1},
(b_{n+1}-a_{n+1})2^{n+1})$, where $\DD(a,r)$ is the open disk
centered at $a$ and of radius $r$.

Observe that for $j\le n$
$$\left|\frac{z-b_j}{z-a_j}-1\right|=\left|\frac{d_j}{z-a_j}\right|<1/2^j$$
on $D_n$. It follows that
$$\prod_{j=1}^n\frac{z-b_j}{z-a_j}$$
is bounded independently of $n$ on $D_n$.

Let  $z_0\in X$. We will show that $(z_0,-f(z_0))\in (\Gamma_f)^*_{\CC^2}$.
Then $\Gamma_{-f}$ is also in the hull and we are done.
Consider the function $g_n$ defined on $\DD \setminus (\cup_{j=1}^n I_j)$ by
$$g_n(z)=\begin{cases} f_n(z) &\text{ if $\Im z<0$;}\\
        -f_n(z) & \text{ if $\Im z>0$;}\\
        \lim_{y\uparrow 0}f(x+iy)&\text{ if $x\in [-1,1]\setminus \cup I_j$.}
\end{cases}$$

The function $g_n$ is holomorphic. Let $c_n=g_n(z_0)+f_K(z_0)$.
Then $c_n\to 0$ as $n\to \infty$. The functions $h_n=g_n-c_n$ tend
to $f_K$ uniformly on compact sets in $\DD\cap \{\Im z<0\}$ and
$h_n(z_0)=-f_K(z_0)$.

Now let $u$ be a plurisubharmonic function on $\CC^2$ that equals $-\infty$ on $\Gamma_f$.
The function
$u(z, h_n(z))$ is subharmonic on the domain $D_n$ and  because $h_n(z)$ is  bounded independently of $n$ on $D_n$, $u(z,h_n(z))$ is bounded by a constant $M$ independently of $n$.

Next we apply the two constant theorem and find
\begin{multline}
u(z_0,-f_K(z_0))\le M(1-\omega(z_0,S,D_n)\\
+\max_{z\in S}u(z, h_n(z))\omega(z_0,S,D_n)\to-\infty,\quad \text{ if $n\to \infty$}.
\end{multline}
Hence, $(\Gamma_f)^*_{\CC^2}\supset\Gamma_f\cup\Gamma_{-f}$.
\par To get the equality we will show that the sup-norm
$$\|g-g_n\|_L^{1/n}\to 0,\quad n\to\infty$$ on compacta $L$
outside $K$. For this we write
$$|g-g_n|=|g_n|
\left| \prod_{k=n+1}^\infty\frac{z-a_k}{z-b_k}-1\right|.$$ The
first factor is bounded by a constant $C$ depending on $L$. To
estimate the second factor we let $\delta$ be the distance from
$L$ to $K$ and write the factor as
$$\left|\prod_{k=n+1}^\infty
\left(1+\frac{d_k}{z-b_k}\right)-1\right|\le\prod_{k=n+1}^\infty
\left(1+\frac{d_k}{\delta}\right)-1\le
\exp\left(\sum\limits_{k=n+1}^\infty
\frac{d_k}{\delta}\right)-1.$$ Since $d_k<4^{-kc_k}$ we see that
$$\left| \prod_{k=n+1}^\infty\frac{z-a_k}{z-b_k}-1\right|\le
\exp\left(\frac{4^{-nc_n}}{(1-4^{-c_n})\delta}\right)-1.$$ Hence
$$\|g-g_n\|_L^{1/n}\le 2\left(\frac{C_L}{\delta}\right)^{1/n}4^{-c_n}$$
when $n$ is sufficiently large and
$$\|g-g_n\|_L^{1/n}\to 0,\quad n\to\infty.$$
\par By Lemma \ref{fa} the pluripolar hull
$\Gamma_g^*=\Gamma_g$ in $\CC^2\setminus(K\times\CC)$. Thus for
any points $(z_0,w_0)$, $z_0\in\CC\setminus K$, $w_0\ne g(z_0)$,
there is a function $u\in\psh(\CC^2)$ such that
$u|_{\Gamma_g}=-\infty$ and $u(z_0,w_0)\ne-\infty$. Then the
function $v(z,w)=u(z,w^2)$ is equal to $-\infty$ on
$\Gamma_f\cup\Gamma_{-f}$ and $v(z_0,\pm\sqrt{w_0})\ne-\infty$.
Hence $(\Gamma_f)^*_{\CC^2}\cap(\CC\setminus
K)\times\CC=\Gamma_f\cup\Gamma_{-f}$.
\end{proof}


\bibliographystyle{amsplain}

\end{document}